\documentclass{article}

\usepackage[all,pdf,2cell]{xy}\UseAllTwocells\SilentMatrices
\usepackage{amsmath}
\usepackage{cleveref}
\usepackage{amsthm}
\usepackage{thmtools}
\usepackage{amsfonts}
\usepackage{enumerate}
\usepackage{todonotes}
\usepackage[dvipsnames]{xcolor}
\usepackage{newunicodechar}
\usepackage{mathtools}
\usepackage{txfonts}
\usepackage{tikz-cd}
\usetikzlibrary{decorations.markings,intersections}
\tikzset{%
scalearrow/.style n args={3}{
  decoration={
    markings,
    mark=at position (1-#1)/2*\pgfdecoratedpathlength
      with {\coordinate (#2);},
    mark=at position (1+#1)/2*\pgfdecoratedpathlength
      with {\coordinate (#3);},
    },
  postaction=decorate,
  }
}
\NewDocumentEnvironment{diagram}{ob}{
  \IfNoValueTF{#1}{\xymatrix{#2}}{\xymatrix@#1{#2}}%
}{}

\newtheorem{theorem}{Theorem}[section]
\newtheorem{lemma}[theorem]{Lemma}
\newtheorem{proposition}[theorem]{Proposition}
\newtheorem{corollary}[theorem]{Corollary}
\theoremstyle{definition}
\newtheorem{definition}[theorem]{Definition}
\newtheorem{example}[theorem]{Example}

\newcommand*{\vcbox}[1]{\begingroup
\setbox0=\hbox{#1}\parbox{\wd0}{\box0}\endgroup}

\newcommand{\id}{\mathrm{id}}

\newcommand{\Aut}{\mathrm{Aut}}

\newcommand{\monicpols}{F_\text{mon}[x]}
\newcommand{\newcategory}[1]{\expandafter\newcommand\csname #1\endcsname{\mathbf{#1}}}
\newcategory{Rel}
\newcategory{RelMonLax}
\newcategory{RelMonOplax}
\newcategory{Set}

\begin{document}
\title{Dimensionalities on monoids}
\author{Gejza Jenča}
\maketitle
\begin{abstract}
We study equivalence relations on monoids in $\Rel$ that are lax morphisms.
We call them \emph{dimensionalities}. We prove that the quotient of
a monoid in $\Rel$ by a dimensionality is again a monoid.
We prove that polygroupoids (and hence hypergroups) are
exactly the symmetric $\dag$-Frobenius structures in $\Rel$.
We use this characterization to prove 
that the quotient of a polygroupoid by a dimensionality that preserves
the inverse is a polygroupoid.
\end{abstract}
\section{Introduction}

In \cite{jencova2017onmonoids} we studied monoids in the category of sets and
relations $\Rel$, which we call relational monoids. We have shown that several
mathematical notions (for example, modular lattices, Riesz decomposition property,
congruences on quantum logic etc.) have a very natural definition when viewed from
the $\Rel$ perspective. In the present paper, we continue this line of research.

We study equivalences on monoids in $\Rel$ that, when considered as
endo-relations, are lax morphisms of monoids. This gives us a new notion even for
ordinary monoids in $\Set$. We call such equivalences \emph{dimensionalities}. We
show that there are naturally occurring  examples, including (but not limited to), equivalences
arising from a subgroup of an automorphism group. The most important property
of dimensionalities is that the quotient of a
relational monoid by a dimensionality is again a relational monoid.

It was proved in \cite{heunen2013relative} that small groupoids are exactly special
$\dag$-Frobenius structures in $\Rel$.  In the last section we extend this result. We
prove that polygroupoids \cite{comer1983new} are exactly symmetric $\dag$-Frobenius structures in $\Rel$.
We then use this result to prove that if a dimensionality preserves
the inverses on a polygroupoid, the quotient by it is again a polygroupoid. As a
consequence, every quotient of a group by a dimensionality that preserves the inverse is a hypergroup.

\section{Preliminaries}

\subsection{The category of sets and relations $\Rel$}

In this section, we review some elementary facts concerning the 2-category of sets
and relations. Everything in this section is well-known, see for example \cite{benabou1967introduction}.

The {\em category of sets and relations}, denoted by $\Rel$, is a category 
whose objects (or $0$-cells) are sets and arrows (or $1$-cells) are relations 
$f\subseteq A\times B$. We will identify such a relation $f$ with a (possibly infinite)
$\{\mathrm{true},\mathrm{false}\}$-valued matrix with columns
indexed by $A$ and rows indexed by $B$, while assuming that
$[\mathrm{true}>\mathrm{false}]$ is a 2-element chain. The fact that $(a,b)\in f$ will be denoted by
$f_a^b$, or by $f(a,b)$.

The composite of arrows $f\colon A\to B$ and 
$g\colon B\to C$ is the arrow $(g\circ f)\colon A\to C$ given by the rule
\begin{equation}\label{eq:composeRel}
(g\circ f)_a^c\Leftrightarrow (\exists b\in B):g_b^c\text{ and }f_a^b.
\end{equation}
Note that the rule \eqref{eq:composeRel} 
may be understood as an analogy of a matrix multiplication rule 
\[
(g\circ f)_a^c=\bigvee_{b\in B}g_b^c\land f_a^b.
\]

The identity arrow $\id_A\colon A\to A$ is the identity relation $\id_A=\{(a,a)\colon a\in A\}$.

There is an obvious faithful functor 
$U\colon\Set\to\Rel$ that is identity on objects and takes a mapping $f\colon A\to B$ to
its graph
$$
U(f)=\{(x,y)\in A\times B\colon f(x)=y\}.
$$
This forgetful functor is a left adjoint, the corresponding right adjoint is the powerset/image functor
$P\colon \Rel\to\Set$. This adjunction induces the well-known {\em powerset monad} on
$\Set$.  $\Rel$ is then isomorphic to the Kleisli category for this monad.

The category of sets and relations is a 2-category: if $h_1,h_2$ are relations $A\to B$, then a
2-cell $h_1\to h_2$ is simply the fact that $h_1\subseteq h_2$. Thus, every hom-category
in $\Rel$ is a poset.

As usual, we draw a 2-cell in a commutative diagram as a double arrow, for example 
\begin{equation}\label{cd:example}
\begin{diagram}
A
	\ar[r]^{f_2}
	\ar[d]_{f_1}
    \ar@{}[rd]|{\Nearrow}
&
B
	\ar[d]^{g_2}
\\
C
	\ar[r]_{g_1}
&
D
\end{diagram}
\end{equation}

means that $g_1\circ f_1\subseteq g_2\circ f_2$. Note that on the level of
elements, the diagram (\ref{cd:example}) means that
\begin{itemize}
\item for every $a\in A$ and $d\in D$ such that there is a $c\in C$ with $(f_1)_a^c$ 
and $(g_1)_c^d$,
\item there exists $b\in B$ such that $(f_2)_a^b$ and $(g_2)_b^d$.
\end{itemize}

Besides the structure of a 2-category, $\Rel$ carries the structure of a 
{\em dagger category}: there is an involution functor $\dag\colon \Rel\to\Rel^{op}$ that is identity on objects.
For a relation $f\colon A\to B$, $f^{\dag}\colon B\to A$ is the relation
given by the rule
\[
(f^{\dag})_b^a=(f)_a^b.
\]

If $A,B\in \Rel$, then the disjoint union of sets $A \sqcup B$ is both the product and the 
coproduct of $A,B$ in $\Rel$. Since $\Rel$ lacks some (co)equalizers,  it is not a (co)complete category.

Note that the dagger structure and the 2-structure on $\Rel$ allow us to express
standard properties of relations on a set: a relation $f\colon A\to A$ is
\begin{itemize}
\item \emph{reflexive}, if $\id\subseteq f$;
\item \emph{symmetric}, if $f=f^\dag$;
\item \emph{transitive}, if $f\circ f\subseteq f$.
\end{itemize}

A relation $\sim$ is an \emph{equivalence}, if it is reflexive, symmetric and
transitive. The quotient of a set $A$ by an equivalence $\sim$ is denoted by $A_\sim$.
For every equivalence $\sim$ there is a quotient mapping $q\colon A\to A_\sim$, given
by $q(a)=[a]_\sim$, where $[a]_\sim$ is the equivalence class of $a$.
Diagrammatically, $q$ can be characterized as the unique mapping such that
\begin{equation}\label{diag:qAndSim}
\xymatrix{
A
    \ar[r]^q
    \ar[rd]_\sim
&
A_\sim
    \ar[rd]^\id
    \ar[d]^{q^\dag}
\\
~
&
A
    \ar[r]_q
&
A_\sim
}
\end{equation}
commutes in $\Rel$.

It is easy to see that the cartesian product $\times$ of sets is a bifunctor from
$\Rel\times\Rel$ to $\Rel$. 

As $\times$ is the product in $\Set$, it satisfies
the coherence conditions for a monoidal category \cite[Chapter VII]{mac1998categories}, 
so $(\Rel,\times,I)$ is a symmetric monoidal category, where $I=\{*\}$ is the
chosen unit element (singleton set) of the monoidal structure.
A \emph{pointed magma} in $\Rel$ is then
an object $S$ equipped with a binary and a nullary operation
\[
\nabla\colon S\times S\to S\qquad e\colon I\to S
\]
A pointed magma in $\Rel$ $(S,\nabla,e)$ is a \emph{monoid} if the diagrams
\begin{center}
\vcbox{
\(
\begin{diagram}
  S\times 1 \ar[r]^{\id\times e} \ar[rd]_{\rho} & S\times S \ar[d]^{\nabla} & 1\times S \ar[l]_{e\times\id} \ar[ld]^{\lambda} \\
  & S &
\end{diagram}
\)
}
$\quad$
\vcbox{
\(
\begin{diagram}
&
S\times S\times S 
    \ar[ld]_{\id\times \nabla}
    \ar[rd]^{\nabla\times\id} \\
  S\times S \ar[rd]_{\nabla} & & S\times S \ar[ld]^{\nabla} \\
  & S &
\end{diagram}
\)
}
\end{center}
commute. The two triangle diagrams are called right and left unit laws. The square
diagram is called the associativity law. Let us note that here and in the rest of
this paper, we omit the associators to avoid unnecessarily complicated diagrams and
computations.

Let us spell out the axioms of a relational monoid in detail. Let $(S,\nabla,e)$ be
a relational monoid. Since $e\colon I\to S$ is a relation, we may identify $e$ with a
subset 
\[
E_S=\{y\in S\colon e_*^y\}
\]
of $S$, which we call {\em the set of units of $S$}.

The $\nabla$ is a relation from $S\times S$ to $S$, so it is a 
subset of $(S\times S)\times S$.
We shall write $\nabla_{ab}^c$ to denote the fact that 
$((a,b),c)\in\nabla$.

The right unit law means that, for every $a\in S$, there is
$y\in E_S$ such that $\nabla_{ay}^a$ and, at the same time,
whenever there is a $y\in E_S$ such that $\nabla_{ay}^b$, then $a=b$. The meaning of
the left unit law is similar.

Associativity law means that for every quadruple $a,b,c,d$ of elements of $S$, the
following statements are equivalent:
\begin{itemize}
\item there exists $p\in S$ such that $\nabla_{ab}^p$ and $\nabla_{pc}^d$;
\item there exists $q\in S$ such that $\nabla_{bc}^q$ and
$\nabla_{aq}^d$.
\end{itemize}

\section{Dimensionalities on relational monoids}

\begin{definition}
Let $S,T$ be pointed magmas in $\Rel$. We say that a relation $h\colon S\to T$ is
a \emph{lax morphism} if
\[
\begin{array}{ccc}
\xymatrix{
S\times S
    \ar[r]^{h\times h}
    \ar[d]_{\nabla}
&
T\times T
    \ar[d]^{\nabla}
\\
S
    \ar[r]^h
    \ar@{}[ru]|{\Nearrow}
    &
T
}
&
\qquad
&
\xymatrix{
I
    \ar[rd]^e
    \ar[d]_e
&
~
\\
S
    \ar[r]_h
    \ar@{}[ru]|(0.35){\Nearrow}
&
T
}\\
(h\circ \nabla)\subseteq (\nabla\circ (h\times h))
&
&
h\circ e\subseteq e
\end{array}
\]
Dually, $h\colon S\to T$ is an \emph{oplax morphism} if and only if
\[
\begin{array}{ccc}
\xymatrix{
S\times S
    \ar[r]^{h\times h}
    \ar[d]_{\nabla}
&
T\times T
    \ar[d]^{\nabla}
\\
S
    \ar[r]^h
    \ar@{}[ru]|{\Swarrow}
    &
T
}
&
\qquad
&
\xymatrix{
I
    \ar[rd]^e
    \ar[d]_e
&
~
\\
S
    \ar[r]_h
    \ar@{}[ru]|(0.35){\Swarrow}
&
T
}\\
(h\circ \nabla)\supseteq (\nabla\circ (h\times h))
&
&
h\circ e\supseteq e
\end{array}
\]
\end{definition}
Let us point out a couple of facts.
\begin{itemize}
\item If $S,T$ are monoids in $\Set$ and $h$ is a mapping, then
$h$ is oplax iff $h$ is lax iff $h$ is a homomorphism of monoids.
\item Let $S$ be a monoid in $\Set$ and let $\sim\colon S\to S$ be an equivalence
relation on $S$. Then $\sim$ is an oplax morphism if and only if $\sim$ is a
congruence.
\item Let $S$ be a monoid in $\Set$ and let $C\subseteq S$. Write $\in_C\colon S\to S$ for the
relation
\[
\in_C(a_1,a_2)\iff a_1=a_2\in C
\]
Then $\in_C$ is an oplax morphism iff $C$ is a submonoid of $S$.
\end{itemize}

We see that several important notions in algebra of monoids can be considered as
particular types of oplax morphisms. General oplax morphisms of semigroups
have been considered in semigroup theory, see the book \cite{rhodes2009theqtheory}.
Nevertheless, in the present paper we will consider lax morphisms, focusing on the
particular case of lax equivalences.

\begin{definition}
Let $(S,\nabla,e)$ be a monoid in $\Rel$. We say that an equivalence relation $\sim$ on
$S$ is a \emph{dimensionality} iff $\sim\colon S\to S$ is a lax morphism.
\end{definition}

This gives us a new notion even in the case of an ordinary monoid in $\Set$.

\begin{proposition}
Let $(A,.,e)$ be a monoid in $\Set$. An equivalence on $A$ is
a dimensionality if the following conditions are satisfied.
\begin{enumerate}[(D1)]
\item If $e\sim e'$, then $e=e'$.
\item For all $a,b,c\in A$ such that $a.b\sim c$, there exist $a',b'\in A$ such that
$a\sim a'$, $b\sim b'$ and $a'.b'=c$.
\end{enumerate}
\end{proposition}
\begin{proof}
Trivial, omitted.
\end{proof}

One large class of examples of dimensionalities comes from automorphism groups acting
on a monoid.
\begin{theorem}\label{thm:automorphisms}
Let $(A,.,e)$ be a monoid in $\Set$, let $\Gamma$ be a subgroup of the group of
automorphisms $\Aut(A)$. Write
\[
x\sim_\Gamma y\text{ iff }\exists f\in\Gamma:f(x)=y
\]
Then $\sim_\Gamma$ is a dimensionality on $A$.
\end{theorem}
\begin{proof}
For (D1), observe that $f(e)=e$ for every $f\in\Aut(A)$. The proof of (D2) is simple,
as well: for every $f\in\Aut(A)$,
\[
f(a.b)=c\implies f(a).f(b)=c
\]
and we may put $a'=f(a)$, $b'=f(b)$.
\end{proof}

Nevertheless, not every dimensionality on a monoid arises in this way, as the following counterexample
shows.
\begin{example}
Let $\Sigma$ be a set with at least 2 elements, consider the monoid 
$(\Sigma^*,.,\epsilon)$ of all words over alphabet $\Sigma$, equipped
with concatenation and the empty word. Let $\sim$ be the 
equivalence on $\Sigma^*$ given by the rule $\mathbf u\sim\mathbf v$ iff
$\mathbf u$ and $\mathbf v$ have the same length. Clearly, $\sim$ is a
dimensionality. However, $\sim$ does not arise from a subgroup of an automorphism
group: we have $aa\sim ab$, yet it is clear that $f(aa)=f(a).f(a)\neq ab$, for any
automorphism $f$.
\end{example}

A natural question arises, whether dimensionalities occur ``in natura''. Of course,
this is partly answered by \Cref{thm:automorphisms}. However, there are several other
examples.

\begin{example}
Write $\monicpols$ for the set of all monic polynomials of one variable over a field $F$:
\[
\monicpols=\{x^n+a_{n-1}x^{n-1}+\dots+a_0:n\in\mathbb N\text{ and }a_{n-1},\dots,a_0\in F\}
\]
Let $\sim$ be an equivalence on $\monicpols$ given by the rule $\mathbf p\sim \mathbf q$ iff 
$\mathbf p$ and $\mathbf q$ have the same degree.

Then $\sim$ is a dimensionality on the monoid $(\monicpols,.,1)$ iff $F$ is
algebraically closed; see \cite[Example 8]{jencova2017onmonoids}
for (essentially) the proof of this.
\end{example}

\begin{example}
Let $(A,+,0)$ be a commutative monoid in $\Set$ and consider the monoid $A^2$. Let $\equiv$ be an equivalence on
$A^2$ given by $(x_1,x_2)\equiv(y_1,y_2)$ iff $x_1+x_2=y_1+y_2$. Then $\equiv$ is a
dimensionality if and only if
\begin{itemize}
\item $A$ is conical, that means satisfies the implication
\[
x+y=0\implies x=y=0
\]
\item $A$ satisfies the common refinement property, that means that $x_1+x_2=y_1+y_2$
implies that there exist $(z_{ij})_{i,j\in\{1,2\}}$ such that
\begin{align*}
x_1=z_{11}+z_{12}&\quad x_2=z_{21}+z_{22}\\
y_1=z_{11}+z_{21}&\quad y_2=z_{12}+z_{22}
\end{align*}
\end{itemize}
Note that $\equiv$ is a congruence on $A^2$.
\end{example}
Examples of such monoids includes all positive cones of partially ordered abelian
groups with interpolation \cite{Goo:POAGwI}. However, there are other meaningful
examples, as well: for example a join-semilattice with zero $(A,\vee,0)$ such
that $\equiv$ is a dimensionality is called a \emph{distributive semilattice}.

Let us explain the reason why we used the name \emph{dimensionality}. For here undefined
notions,
we refer to the book \cite{DvuPul:NTiQS} and the paper \cite{jencova2017onmonoids}.
\begin{definition}
\cite{loomis1955lattice}
Let $A$ be a complete orthomodular lattice. A {\em dimension equivalence} on $A$
is an equivalence relation on $A$ such that
\begin{enumerate}
\item[(A)]If $a\sim 0$, then $a=0$.
\item[(B)]If $a_1\perp a_2$ and $a_1\vee a_2\sim b$, then there exists an orthogonal
decomposition of $b$, $b=b_1\vee b_2$, such that $b_1\sim a_1$ and $b_2\sim a_2$.
\item[(C)]If $\{a_\alpha\}$ and $\{b_\alpha\}$ are pairwise orthogonal families of elements,
such that $a_\alpha\sim b_\alpha$ for all $\alpha$, then $\bigvee_\alpha a_\alpha=\bigvee_\alpha b_\alpha$.
\item[(D)]If $a$ and $b$ are not orthogonal in $A$ then there are nonzero $a_1,b_1$ in $A$ such that
$a\geq a_1$, $b\geq b_1$ and $a_1\sim b_1$.
\end{enumerate}
\end{definition}

Every orthomodular lattice is a partial abelian monoid, if we restrict the operation
$\vee$ to orthogonal pairs of elements. For this partial abelian monoid, the 
conditions (A) and (B) in the previous definition exactly mean that
a dimension equivalence is a dimensionality.

\section{Quotients}

For every equivalence $\sim$ on a pointed magma $(S,\nabla,e)$ in $\Rel$ we may
construct a pointed magma $(S_\sim,\nabla_\sim,e_\sim)$, where $\nabla_\sim$ and
$e_\sim$ are given by the rules 
\begin{equation}\label{eq:simOps}
\begin{split}
(\nabla_\sim)_{AB}^C & \iff \exists(a,b)\in A\times B,c\in C:\nabla_{ab}^c\\
(e_\sim)_*^P&\iff \exists p\in P:e_*^p
\end{split}
\end{equation}

Then \eqref{diag:qAndSim} gives us a canonical
factorization $\sim=q^\dag\circ q$ and the formulas in \eqref{eq:simOps} easily imply
that
\begin{equation}\label{diag:qIsLax}
\xymatrix{
S\times S
    \ar[r]^-{q\times q}
    \ar[d]_{\nabla}
&
S_\sim\times S_\sim
    \ar[d]^{\nabla_\sim}
\\
S
    \ar[r]_q
    \ar@{}[ru]|{\Nearrow}
&
S_\sim
}
\qquad
\xymatrix{
&
1
    \ar[ld]_e
    \ar[rd]^{e_\sim}
    \ar@{}[d]|(0.65){\Nearrow}
\\
S
    \ar[rr]_q
&
&
S_\sim
}
\end{equation}
or, in other words, that $q$ is a lax morphism of pointed magmas.

\begin{proposition}\label{prop:dimenIffQDagIsLax}
An equivalence on a relational monoid is a dimensionality if and only if
$q^\dag$ is a lax morphism of pointed magmas.
\end{proposition}
\begin{proof}
Pasting the diagrams \eqref{diag:qIsLax} with the analogous
diagrams for $q^\dag$ gives us
\begin{equation}\label{diag:qAndQDagAreLax}
\xymatrix{
S\times S
    \ar[r]^{q\times q}
    \ar[d]_{\nabla}
&
S_\sim\times S_\sim
    \ar[r]^{q^\dag\times q^\dag}
    \ar[d]_{\nabla_\sim}
&
S\times S
    \ar[d]^\nabla
\\
S
    \ar[r]_q
    \ar@{}[ru]|{\Nearrow}
&
S_\sim
    \ar[r]_{q^\dag}
    \ar@{}[ru]|{\Nearrow}
&
S
}
\qquad
\xymatrix{
&
1
    \ar[ld]_e
    \ar[rd]^{e_\sim}
    \ar[rr]^e
    \ar@{}[d]|(0.65){\Nearrow}
&
&
S
\\
S
    \ar[rr]_q
&
&
S_\sim
    \ar[ru]_{q^\dag}
    \ar@{}[u]|(0.65){\Nearrow}
}
\end{equation}
Suppose that $q^\dag$ is a lax morphism. Then
all of the inner diagrams in \eqref{diag:qAndQDagAreLax} commute. Clearly, this
implies that $\sim=q^\dag\circ q$ is a lax morphism and hence a dimensionality.

On the other hand, suppose that $\sim=q^\dag\circ q$ is a dimensionality.  The
commutativity of the $q^\dag$ square in \eqref{diag:qAndQDagAreLax} exactly means
that for all $A,B\in S_\sim$ and $c'\in S$ such that there exists $a\in A$, $b\in B$
and $c\sim c'$ such that $\nabla_{ab}^c$, there exist $a'\in A$, $b'\in B$ such that
$\nabla_{a'b'}^{c'}$, but this is exactly
$(\sim\circ\nabla_\sim)\subseteq\bigl(\nabla\circ(\sim\times\sim)\bigr)$. 
Similarly, the right-hand side triangle in \eqref{diag:qAndQDagAreLax} exactly means that $(\sim\circ e_\sim)\subseteq e$.
\end{proof}

\begin{theorem}\label{thm:dimensionality}
For every dimensionality on a relational monoid $(S,\nabla,e)$, 
$(S_\sim,\nabla_\sim,e_\sim)$ is a relational monoid.
\end{theorem}
\begin{proof}
Consider the diagram
\begin{equation}\label{diag:transferAssoc}
\xymatrix{
&
(S_\sim)^3
    \ar[ld]_-{\id\times\nabla_\sim}
    \ar[r]^-{(q^\dag)^3}
    \ar@/^2em/[rr]^\id
    \ar@{}[d]|{\Nearrow}
&
S^3
    \ar[ld]^-{\id\times\nabla}
    \ar[r]^-{q^3}
    \ar[rd]_-{\nabla\times\id}
&
(S_\sim)^3
    \ar[rd]^-{\nabla_\sim\times\id}
    \ar@{}[d]|{\Nearrow}
\\
(S_\sim)^2
    \ar[rd]_{\nabla_\sim}
    \ar[r]_-{(q^\dag)^2}
&
S^2
    \ar[rd]^\nabla
    \ar@{}[d]|{\Nearrow}
&
&
S^2
    \ar[ld]_\nabla
    \ar[r]_-{q^2}
    \ar@{}[d]|{\Nearrow}
&
(S_\sim)^2
    \ar[ld]^{\nabla_\sim}
\\
&
S_\sim
    \ar[r]_{q^\dag}
    \ar@/_2em/[rr]_-{\id}
&
S
    \ar[r]_q
&
S_\sim
}
\end{equation}
By \eqref{diag:qAndSim}, the top and bottom triangles commute. The left-hand side lax
squares commute by \Cref{prop:dimenIffQDagIsLax}, the right-hand side lax squares
commute by \eqref{diag:qIsLax}. The middle square is the associativity law for
$\nabla$. The outside border of the diagram gives us
\[
\nabla_\sim\circ(\id\times\nabla_\sim)\subseteq(\nabla_\sim\times\id)\circ\nabla_\sim
\]
The reverse inclusion follows in a similar way, all that is needed is to flip the
components in the direct product in the top half of \eqref{diag:transferAssoc}.

The two inclusions needed for the left unit law follow by
\begin{equation}\label{diag:transferUnit1}
\xymatrix{
&
I\times S_\sim
    \ar[r]^-{\id\times q^\dag}
    \ar@/^2em/[rr]^-{\id}
    \ar[ld]_{e_\sim\times\id}
    \ar@{}[d]|{\Nearrow}
&
I\times S
    \ar[r]^-{\id\times q}
    \ar[ld]^-{e\times\id}
    \ar[dd]^{\lambda}
    %\ar@{}[rdd]|{\Nearrow}
&
I\times S_\sim
    \ar[dd]^{\lambda}
\\
(S_\sim)^2
    \ar[r]_{q^\dag\times q^\dag}
    \ar[rd]_-{\nabla_\sim}
&
S^2
    \ar[rd]^-\nabla
    \ar@{}[d]|{\Nearrow}
\\
&
S_\sim
    \ar[r]_-{q^\dag}
    \ar@/_2em/[rr]_-{\id}
&
S
    \ar[r]_-{q}
&
S_\sim
}
\end{equation}
\begin{equation}\label{diag:transferUnit2}
\xymatrix{
I\times S_\sim
    \ar[r]^{\id\times q^\dag}
    \ar@/^2em/[rr]^-{\id}
    \ar[dd]_{\lambda}
&
I\times S
    \ar[r]^{\id\times q}
    \ar[rd]_-{e\times\id}
    \ar[dd]_{\lambda}
&
I\times S_\sim
    \ar[rd]^-{e_\sim\times\id}
    \ar@{}[d]|{\Nearrow}
\\
&
&
S^2
    \ar[r]_{q\times q}
    \ar[ld]_-\nabla
    \ar@{}[d]|{\Nearrow}
&
(S_\sim)^2
    \ar[ld]^-{\nabla_\sim}
\\
S_\sim
    \ar@/_2em/[rr]_-{\id}
    \ar[r]_{q^\dag}
&
S
    \ar[r]_{q}
&
S_\sim
}
\end{equation}
Here, $\lambda$ denotes the left unitor of the monoidal structure.
In \eqref{diag:transferUnit1}, the top left square is the product of the lax square 
\[
q^\dag\circ e_\sim\subseteq e=e\circ\id
\]
where the inclusion follows by \Cref{prop:dimenIffQDagIsLax},
and the trivial square $q^\dag\circ\id=\id\circ q^\dag$. The bottom lax square
follows by \Cref{prop:dimenIffQDagIsLax} as well. The triangle in the middle
is the left unit law for $\nabla$ and $e$. The strict square is the naturality of
$\lambda$. The top and bottom triangles are from \eqref{diag:qAndSim}.

The proof that the cells in \eqref{diag:transferUnit2} commute is very similar. These diagrams,
when taken together, give us the left unit law. The proof of the right unit law follows by a
symmetrical argument.
\end{proof}

Utilizing \Cref{thm:dimensionality} and \Cref{thm:automorphisms}, we can easily
construct examples of relational monoids from ordinary monoids in $\Set$ and
subgroups of their automorphism groups, as in the following example.

\begin{example}\label{ex:rotations}
Consider the commutative monoid of all vectors $(\mathbb R^2,+,(0,0))$ and the group
$\Gamma$ of all $2\times 2$ rotation matrices, acting on $\mathbb R^2$ in the usual
way. Then two vectors are related by $\sim_\Gamma$ iff they have the same length,
so we may identify $\mathbb R^2/\sim_\Gamma$ with the nonnegative reals $[0,\infty)$ and 
then the relation $+\colon [0,\infty)\times[0,\infty)\to[0,\infty)$ is given by the rule
\[
+((a,b),x)\iff |a-b|\leq x\leq a+b
\]
\end{example}

Let us extend \Cref{thm:automorphisms} to relational monoids. The first question we
have to ask is about the correct notion of isomorphism of relational
monoids. Clearly, an isomorphism should be an invertible lax (or oplax) morphism.
Invertible relations coincide with bijective maps, so the notion we are seeking must be a lax or oplax bijective morphism such that the inverse is again a lax or oplax
morphism. The following lemma and its corollary show that
these are, actually, the same.
\begin{lemma}\label{lemma:OplaxLaxIso}
Let $S$, $T$ be relational monoids, let $h\colon S\to T$ be a bijection.
Then $h$ is a lax morphism of monoids if and only if $h^{-1}$ is an oplax morphism of
monoids.
\end{lemma}
\begin{proof}
It is easy to see that 
$h$ is lax if and only if, for all $y,s_1,s_2\in S$ and 
$t\in T$, 
\begin{align}\label{eq:isoLaxNablaMult}
\nabla_{s_1,s_2}^{h^{-1}(t)}\implies\nabla_{h(s_1),h(s_2)}^t\\
\label{eq:isoLaxNablaUnit}
y\in E_S\implies h(y)\in E_T
\end{align}
Similarly, $h^{-1}$ is oplax if and only if for all $t_1,t_2\in T$, $y,s\in S$,
\begin{align}\label{eq:isoOplaxNablaMult}
\nabla_{h^{-1}(t_1),h^{-1}(t_2)}^{s}\implies\nabla_{t_1,t_2}^{h(s)}\\
\label{eq:isoOplaxNablaUnit}
y\in E_S\implies\exists w\in E_T:h^{-1}(w)=y
\end{align}
Assume that $h$ is lax. 
Put $s_1:=h^{-1}(t_1)$, $s_2:=h^{-1}(t_2)$ and $t:=h(s)$ in \eqref{eq:isoLaxNablaMult}
to obtain \eqref{eq:isoOplaxNablaMult}. Put $w:=h(y)$ in \eqref{eq:isoLaxNablaUnit}
to obtain \eqref{eq:isoOplaxNablaUnit}.

Assume that $h^{-1}$ is oplax. 
Put $t_1:=h(s_1)$, $t_2:=h(s_2)$ and $s:=h^{-1}(t)$ in \eqref{eq:isoOplaxNablaMult}
to obtain \eqref{eq:isoLaxNablaMult}. 
To prove \eqref{eq:isoLaxNablaUnit}, let $y\in E_S$. By \eqref{eq:isoOplaxNablaUnit},
there is $w\in E_T$ with $h^{-1}(w)=y$. But if $h^{-1}(w)=y$, then $h(y)=w\in E_T$.
\end{proof}
\begin{corollary}\label{coro:autoLaxIsOplax}
Let $S$, $T$ be relational monoids, let $h\colon S\to T$ be a bijection.
Then both $h$ and $h^{-1}$ are lax if and only if both $h$ and $h^{-1}$ are oplax.
\end{corollary}
\begin{proof}
Apply \Cref{lemma:OplaxLaxIso} twice.
\end{proof}
Hence the isomorphisms (automorphisms) in the categories lax and
oplax morphisms of relational monoids coincide, and we may simply call them \emph{isomorphisms
(automorphisms)}, without a lax/oplax qualifier.

\begin{theorem}
Let $(S,\nabla,e)$ be a relational monoid, let $\Gamma$ be a subgroup of the group of
automorphisms $\Aut(S)$. Write
\[
x\sim_\Gamma y\text{ iff }\exists f\in\Gamma:f(x)=y
\]
Then $\sim_\Gamma$ is a dimensionality.
\end{theorem}
\begin{proof}
We need to prove that $\nabla_{ab}^c$ and $c\sim_\Gamma c'$ imply that
there are $a',b'\in S$ such that $\nabla_{a'b'}^{c'}$. There is $h\in\Gamma$
such that $c'=h(c)$. Since $h$ is an isomorphism, $h$ is oplax, hence $\nabla_{ab}^c$
implies $\nabla_{h(a)h(b)}^{h(c)}$ and we may put $a'=h(a)$ and $b'=h(b)$. The
$\sim_\Gamma\circ e\subseteq e$ condition follows in a similar way.
\end{proof}

\section{Frobenius structures, polygroupoids and hypergroups}

A \emph{Frobenius structure} in  $\Rel$ is
a structure $(S,\nabla,\Delta,e,c)$ such that $(S,\nabla,e)$ is a pointed magma,
$\Delta\colon S\to S\times S$, $c\colon S\to I$ are relations and the following diagrams commute.
\begin{center}
\vcbox{
\(
\begin{diagram}
  S\times I \ar[r]^-{\id\times e} \ar[rd]_{\rho} & S\times S \ar[d]^{\nabla} & I\times
  S \ar[l]_-{e\times\id} \ar[ld]^{\lambda} \\
  & S &
\end{diagram}
\)
}
\vcbox{
\(
\begin{diagram}
  & S
  \ar[ld]_{\rho^{-1}}
  \ar[d]^\Delta
  \ar[rd]^{\lambda^{-1}}
  \\
  S\times I 
  &
  S\times S
  \ar[l]^-{\id\times c}
  \ar[r]_{c\times\id}
  &
  I\times S
\end{diagram}
\)
}
\end{center}
\begin{equation}
\label{diag:frobenius}
\xymatrix{
S\times S
	\ar[rr]^-{\id\times\Delta}
	\ar[rd]^-{\nabla}
	\ar[dd]_-{\Delta\times\id}
&
~
&
S\times S\times S
	\ar[dd]^-{\nabla\times\id}
\\
~
&
S
	\ar[rd]^-{\Delta}
\\
S\times S\times S
	\ar[rr]_-{\id\times\nabla}
&
~
&
S\times S
}
\end{equation}

The top left diagram are just the familiar unit laws. The top right diagram are the 
\emph{counit laws}. The bottom diagram is the \emph{Frobenius law}.
By \cite[Lemma 3.12]{mehta2020frobenius}, $\nabla$ and $\Delta$ in a Frobenius structure
are (co)associative. Hence every Frobenius structure is a monoid and a comonoid.

A Frobenius structure is \emph{symmetric} if and only if
$c\circ\nabla=c\circ\nabla\circ\sigma$.

By \cite[Proposition 3.4]{mehta2020frobenius}, there is a bijection $\hat\beta\colon S\to S$
that is completely determined by the property $\nabla_{a\hat\beta(a)}^p$ for some $p\in
E_S$. Moreover, by \cite[Lemma 3.7]{mehta2020frobenius}
\begin{equation}\label{eq:betaNoDag}
\nabla_{xy}^z\iff\Delta_y^{\hat\beta(x)z}\iff\Delta_x^{z({\hat\beta})^{-1}(y)}.
\end{equation}

We say that a Frobenius structure is a \emph{$\dag$-Frobenius structure} if $c=e^\dag$ and
$\Delta=\nabla^\dag$. For a $\dag$-Frobenius structure, \eqref{eq:betaNoDag} is
equivalent to
\[
\nabla_{xy}^z\iff\nabla^y_{\hat\beta(x)z}\iff\nabla^x_{z{\hat\beta}^{-1}(y)}.
\]

\begin{lemma}\label{lemma:nondegBijection}
A $\dag$-Frobenius structure in $\Rel$ is symmetric if and only if
$\hat\beta\circ\hat\beta=\id$.
\end{lemma}
\begin{proof}
Recall that $\hat\beta$ is determined by the property that, for all $x,y\in S$,
\begin{equation}\label{eq:betaDef}
(\exists p\in E_S:\nabla_{xy}^p)\iff y=\hat\beta(x).
\end{equation}
Since $c=e^\dag$, the structure is symmetric if and only if
$e^\dag\circ\nabla=e^\dag\circ\nabla\circ\sigma_{S,S}$, that is, if and only if for all
$x,y\in S$
\[
(\exists p\in E_S:\nabla_{xy}^p)\iff(\exists q\in E_S:\nabla_{yx}^q).
\]
By \eqref{eq:betaDef}, the left-hand side is $y=\hat\beta(x)$ and the right-hand side
is $x=\hat\beta(y)$, so symmetry is equivalent to
\begin{equation}\label{eq:betaInvol}
y=\hat\beta(x)\iff x=\hat\beta(y)\qquad\text{for all }x,y\in S.
\end{equation}
Since $\hat\beta$ is a bijection, $\hat\beta\circ\hat\beta=\id$ if and only if
$\hat\beta=\hat\beta^{-1}$; and as $x=\hat\beta(y)$ is equivalent to
$y=\hat\beta^{-1}(x)$, the identity $\hat\beta=\hat\beta^{-1}$ is precisely
\eqref{eq:betaInvol}. Hence the structure is symmetric if and only if
$\hat\beta\circ\hat\beta=\id$.
\end{proof}
The following definition is just a recasting of the original definition from \cite{comer1983new} to
the language of $\Rel$.
\begin{definition}\label{def:polygroupoid}
A \emph{polygroupoid} is a relational monoid $(S,\nabla,e)$ equipped with a mapping
$x\mapsto x^{-1}$ such that 
\begin{equation}\label{eq:nablaFlip}
\nabla_{xy}^z\iff \nabla_{zy^{-1}}^x\iff\nabla_{x^{-1}z}^y
\end{equation}
\end{definition}
A polygroupoid is a \emph{groupoid} iff $\nabla$ is a partial operation. A polygroupoid
is a \emph{hypergroup} iff $E_S$ is a singleton. 

The following is just a rewriting parts of \cite[Lemma 4.1]{comer1983new} to relational
language.
\begin{lemma}\label{lemma:comer}
Let $(S,\nabla,e,~^{-1})$ be a polygroupoid. Then the following hold for all
$a,b\in S$:
\begin{enumerate}[(a)]
\item There are $p,q\in E_S$ such that $\nabla_{aa^{-1}}^p$ and $\nabla_{a^{-1}a}^q$.
\item If $\nabla_{ab}^p$ for some $p\in E_S$, then $b=a^{-1}$.
\item $(a^{-1})^{-1}=a$.
\end{enumerate}
\end{lemma}

\begin{corollary}\label{coro:autoPreservesInverse}
For every monoid automorphism $f$ of a polygroupoid $S$, $f(a^{-1})=(f(a))^{-1}$
for all $a\in S$.
\end{corollary}
\begin{proof}
For every $a\in S$, there is a $p\in E_S$ such that $\nabla_{aa^{-1}}^p$ by (a).
Since $f$ is an automorphism, $f$ is oplax by \Cref{coro:autoLaxIsOplax}. Therefore,
$\nabla_{f(a)f(a^{-1})}^{f(p)}$ and $f(p)\in E_S$. By (b) and (c), we see that
$(f(a))^{-1}=f(a^{-1})$.

\end{proof}
%We say that a 2-form $\beta\colon S\otimes S\to I$ 
%is nodegenerate if and only if it is a part of some self-duality $(S,S,\eta,\beta)$.
%By \cite[Lemma 3.2]{mehta2020frobenius}, a 2-form
%$\beta\colon S\otimes S\to I$ in $\Rel$ is nondegenerate if and only if there is
%bijective map $\hat\beta\colon S\to S$ such that
%\[
%\beta_{xy}^*\iff y=\hat\beta(x)
%\]

%We say that a 2-form $\beta\colon A\to A$ is symmetric iff
%$\beta=\beta\circ\sigma_{A,A}$.
%\begin{lemma}\label{lemma:nondegBijection}
%A nondegenerate 2-form $\beta\colon X\otimes X\to I$ is symmetric if and only if
%$\hat\beta\circ\hat\beta=\id$.
%\end{lemma}
%\begin{proof}
%$\hat\beta\circ\hat\beta=\id$ 
%is equivalent to $y=\hat\beta(x)$ iff $x=\hat\beta(y)$, but this exactly
%means that $\beta$ is symmetric.
%\end{proof}
\begin{theorem}\label{thm:polygroupoidIsFrob}
There is a one-to-one correspondence between polygroupoids 
and symmetric dagger Frobenius structures in $\Rel$.
\end{theorem}
\begin{proof}
Suppose that $(S,\nabla,e,~^{-1})$ is a polygroupoid in $\Rel$.
Then $(S,\nabla,e)$ is a monoid in $\Rel$, so $(S,\nabla^\dag,e^\dag)$ is
a comonoid in $\Rel$. 
Putting $\Delta:=\nabla^\dag$, \eqref{eq:nablaFlip} becomes
\begin{equation}\label{eq:deltaFlip}
\nabla_{xy}^z\iff\Delta^{zy^{-1}}_x\iff\Delta^{x^{-1}z}_y.
\end{equation}
Let us prove this implies the Frobenius law for $\nabla$ and $\Delta$.
For $(a,b),(c,d)\in S\times S$,
\[
\bigvee_{p\in S}\nabla_{ap}^c\wedge\Delta_b^{pd}=
\bigvee_{p\in S}\nabla_{ap}^c\wedge\nabla_{bd^{-1}}^p=
\bigvee_{q\in S}\nabla_{ab}^q\wedge\nabla_{qd^{-1}}^c=
\bigvee_{q\in S}\nabla_{ab}^q\wedge\Delta_q^{cd}
\]
proving one of the Frobenius laws. The other one follows in an analogous way.
To prove that this Frobenius structure is symmetric, we need to prove that, for all $a,b\in S$,
\[
(e^\dagger\circ\nabla)_{ab}^*=(e^\dagger\circ\nabla\circ\sigma_{S,S})_{ab}^*
\]
This amounts to proving that
\[
(\exists p\in E_S:\nabla_{ab}^p)
\iff
(\exists q\in E_S:\nabla_{ba}^q)
\]
By symmetry, it suffices to prove just one implication, say $\implies$.
This follows by \Cref{lemma:comer}:
by (b), we see that $\nabla_{ab}^p$ and $p\in E_S$ means that $b=a^{-1}$. By (c), this is equivalent
to $a=b^{-1}$ and by (a), this implies $\nabla_{ba}^q=\nabla_{a^{-1}a}^q$ for some $q\in E_S$.

On the other hand, suppose that we have a symmetric dagger Frobenius structure $A$ in
$\Rel$. By \cite[Proposition 3.4]{mehta2020frobenius} there is
a bijection $\hat\beta\colon A\to A$ such that
\[
(\exists p\in E_S:\nabla_{xy}^p)\iff y=\hat\beta(x)
\]
As our Frobenius structure is symmetric, \Cref{lemma:nondegBijection} tells us
that $(\hat\beta)^{-1}=\hat\beta$.
By \cite[Lemma 3.7]{mehta2020frobenius}, $\hat\beta$ satisfies
\[
\nabla_{xy}^z\iff\Delta_y^{\hat\beta(x)z}\iff\Delta_x^{z({\hat\beta})^{-1}(y)}.
\]
Therefore, as $\Delta=\nabla^\dag$,
\[
\nabla_{xy}^z\iff\nabla^y_{\hat\beta(x)z}\iff\nabla^x_{z\hat\beta(y)}.
\]
and we see that \eqref{eq:nablaFlip} is satisfied for $\hat\beta=~^{-1}$.
\end{proof}

\begin{theorem}\label{thm:polygroupoid}
Let $(S,\nabla,e)$ be a polygroupoid, let $\sim$
be a dimensionality on $S$ such that, for all $x,y\in S$,
\begin{equation}\label{eq:invPreserved}
x\sim y\iff x^{-1}\sim y^{-1}
\end{equation}
Then $(S_\sim,\nabla_\sim,e_\sim)$ is a polygroupoid. 
\end{theorem}
\begin{proof}
By \Cref{thm:dimensionality}, the quotient $(S_\sim,\nabla_\sim,e_\sim)$ is a
relational monoid. Since $\sim$ preserves the inversion, the inversion mapping 
can be transferred to $S_\sim$ and satisfies \eqref{eq:nablaFlip}. 
\end{proof}
\begin{corollary}\label{coro:hyperGroupPreserved}
If $(S,\nabla,e)$ is a hypergroup and $\sim$ is a dimensionality on $S$ that
satisfies \eqref{eq:invPreserved}, then $(S_\sim,\nabla_\sim,e_\sim)$ is a
hypergroup.
\end{corollary}
\begin{corollary}
If $(S,\nabla,e)$ is a symmetric dagger Frobenius structure in $\Rel$
and $\sim$ is a dimensionality on $S$ that satisfies \eqref{eq:invPreserved},
then $(S_\sim,\nabla_\sim,e_\sim)$ is a symmetric dagger Frobenius structure.
\end{corollary}
\begin{proof}
Directly from \Cref{thm:polygroupoid} and \Cref{thm:polygroupoidIsFrob}.
\end{proof}
\begin{example}
By \Cref{coro:autoPreservesInverse} and \Cref{coro:hyperGroupPreserved}, the
relational monoid we obtained in \Cref{ex:rotations} is a hypergroup. This hypergroup
forms the additive part of the \emph{triangle hyperfield}, see
\cite{viro2010hyperfields}.
\end{example}
\begin{example}
Consider the symmetric group $S_n$. As every group, $S_n$ acts on itself by a
conjugation, via inner automorphisms. The orbits of this action correspond to cycle
types of the permutations, and \Cref{coro:hyperGroupPreserved} means that 
they do form a hypergroup. For $n=3$ we obtain cycle types $\{(),(2),(3)\}$
and the hypergroup of cycle types is
\[
\begin{array}{c|c|c|c}
    &() &(2)    &(3)    \\
\hline
()  &() &(2)    &(3)    \\
(2) &(2)&();(3) &(2)    \\
(3) &(3)&(2)    &();(3)
\end{array}
\]
For $n=4$ we obtain cycle types $\{(),(2),(3),(2,2),(4)\}$ and the hypergroup 
of cycle types is
\[
\begin{array}{c|c|c|c|c|c}
      & ()     & (2)          & (2,2)      & (3)          & (4)          \\
\hline
()    & ()     & (2)          & (2,2)      & (3)          & (4)          \\
(2)   & (2)    & ();(2,2);(3) & (2);(4)    & (2);(4)      & (2,2);(3)    \\
(2,2) & (2,2)  & (2);(4)      & ();(2,2)   & (3)          & (2);(4)      \\
(3)   & (3)    & (2);(4)      & (3)        & ();(2,2);(3) & (2);(4)      \\
(4)   & (4)    & (2,2);(3)    & (2);(4)    & (2);(4)      & ();(2,2);(3)
\end{array}
\]
Note that both hypergroups are commutative. Indeed, this is true for every
$S_n$. To see this, let us write $f\sim g$ for having the same cycle type,
so that $[f]_\sim$ is a cycle type of $f\in S_n$. 
For every $f\in S_n$ we have $[f]_\sim=[f^{-1}]_\sim$.
Let $F,G,H\in
(S_n)_\sim$ and suppose that $(\nabla_\sim)_{FG}^H$, meaning that 
\[
\exists(f,g)\in F\times G,h\in H:fg=h
\]
Then $g^{-1}f^{-1}=h^{-1}$ and this implies $\nabla_{GF}^H$.
The reverse implication follows by symmetry.
\end{example}
\begin{example}
Let $F$ be a field. Consider a subgroup $G$ of the multiplicative group 
$F^\times=(F\setminus\{0\},.,1)$ of $F$. The group $G$ acts on the additive group
of $F$ by multiplication and preserves the mapping $x\mapsto -x$. Therefore,
$(F/\sim_G,+,0)$ is an (abelian) hypergroup. We remark that much more is true in this
example: $F/\sim_G$ is a \emph{hyperfield}. 
\end{example}
\section*{Declarations}
\begin{description}
\item[\bf Funding:]
This research is supported by grants VEGA 2/0128/24 and 1/0036/23,
Slovakia.
\end{description}

%\bibliography{mybib}

\begin{thebibliography}{10}

\bibitem{benabou1967introduction}
Jean B{\'e}nabou.
\newblock Introduction to bicategories.
\newblock In {\em Reports of the Midwest Category Seminar}, pages 1--77. Springer, 1967.

\bibitem{comer1983new}
Stephen~D Comer.
\newblock A new foundation for the theory of relations.
\newblock {\em Notre Dame Journal of Formal Logic}, 24(2):181--187, 1983.

\bibitem{DvuPul:NTiQS}
A.~Dvure{\v c}enskij and S.~Pulmannov{\'a}.
\newblock {\em New Trends in Quantum Structures}.
\newblock Kluwer, Dordrecht and Ister Science, Bratislava, 2000.

\bibitem{Goo:POAGwI}
K.R. Goodearl.
\newblock {\em Partially ordered abelian groups with interpolation}.
\newblock Amer. Math. Soc, Providence, 1986.

\bibitem{heunen2013relative}
Chris Heunen, Ivan Contreras, and Alberto~S Cattaneo.
\newblock Relative {F}robenius algebras are groupoids.
\newblock {\em Journal of Pure and Applied Algebra}, 217:114--124, 2013.

\bibitem{jencova2017onmonoids}
Anna Jen{\v{c}}ov{\'a} and Gejza Jen{\v{c}}a.
\newblock On monoids in the category of sets and relations.
\newblock {\em International Journal of Theoretical Physics}, 56(12):3757--3769, 2017.

\bibitem{mac1998categories}
Saunders~Mac Lane.
\newblock {\em Categories for the Working Mathematician}.
\newblock Number~5 in Graduate Texts in Mathematics. Springer-Verlag, 1971.

\bibitem{loomis1955lattice}
Lynn~H Loomis.
\newblock The lattice theoretic background of the dimension theory of operator algebras.
\newblock {\em Memoirs of the AMS}, 18, 1955.

\bibitem{mehta2020frobenius}
Rajan~Amit Mehta and Ruoqi Zhang.
\newblock Frobenius objects in the category of relations.
\newblock {\em Letters in Mathematical Physics}, 110(7):1941--1959, 2020.

\bibitem{rhodes2009theqtheory}
John~L Rhodes and Benjamin Steinberg.
\newblock {\em The q-theory of finite semigroups}.
\newblock Springer, 2009.

\bibitem{viro2010hyperfields}
Oleg Viro.
\newblock Hyperfields for tropical geometry {I}. {H}yperfields and dequantization.
\newblock {\em arXiv preprint arXiv:1006.3034}, 2010.

\end{thebibliography}
%\bibliographystyle{plain}
\end{document}